\magnification=\magstep1

\newcount\sec \sec=0
\input Ref.macros
\input math.macros
\forwardreferencetrue
\citationgenerationtrue
\initialeqmacro
\sectionnumberstrue

\def\hg{hypergraph }

\def\mlp{minimal cover set }
\def\mcs{minimal cover set }

\title{Split hypergraphs}
\author{\'Ad\'am Tim\'ar}
\abstract{Generalizing the notion of split graphs to uniform hypergraphs, 
we prove
that the class of these hypergraphs can be characterized by a finite list
of excluded induced subhypergraphs. We show that a characterization by 
generalized degree sequences is impossible, unlike in the well-known case 
of split graphs. We also give an algorithm to decide
whether a given uniform hypergraph is a split hypergraph. If it is, the 
algorithm gives a
splitting of it; the running time is $O(N\log N)$. These answer 
questions of 
Sloan, Gy. Tur\'an and Peled.}
\bottomII{Primary 05C65.  Secondary 05C85, 68R05. }
{Split hypergraph, split graph, excluded induced subgraph, degree 
sequence, sparse and dense graph.}
{Research done when the author was at the University of Szeged and 
partially supported by Hungarian National Foundation for Scientific 
Research Grant OTKA49398.}

\bsection{Introduction, split graphs}{s.intro} 

We call a graph $G$ {\it split} if there exists a partition $A\cup B$ of
its vertex set so that there is no edge between any two points of $A$ and 
there is an edge between every two points of $B$.
Such a partition $A\cup B$ of
$V(G)$ is called a {\it splitting}.  

The following theorem describes split graphs in two different ways:

\procl t.graph
For a graph $G$ the following are equivalent:
 
(1) $G$ is split.

(2) $G$ does not contain $C_4$, $C_5$ or $K_2 \cup K_2$ as an induced 
subgraph.

(3) There is an $m$ such that $\sum^m_{i=1} d_i - m(m-1) = \sum^n_{i=m+1}
d_i$, where $d_1\ge d_2 \ge \ldots \ge d_n$ are the degrees in $G$. 
\endprocl

The equivalence of (1) and (2) is first shown in \ref b.GyL/, while the 
equivalence of (1) and (3) is in \ref b.HS/.
See also \ref b.BLS/ or \ref b.Go/ for history and further 
references.

The generalizations of many graph theoretical questions to hypergraphs
have practical significance. Most frequently, posing these problems for
hypergraphs makes them rather difficult. For example, testing the
existence of a perfect matching can be solved by a well-known
polynomial algorithm for graphs. Generalizing it to 
3-uniform
hypergraphs we get an NP-complete problem (\ref b.K/, \ref b.GJ/).

The class of line-graphs can be characterized by a finite list of
excluded induced subgraphs --- on the other hand there are infinitely
many pairwise nonisomorphic graphs with the property that none of them is 
the line-graph of any
3-uniform hypergraph, but every proper subgraph of any of them is, \ref 
b.NRSS/. 
Studying extremal
questions for graphs turned out to be very fruitful. The generalization of
this theory to hypergraphs comes up against serious difficulties.

Recent research in learning theory has led to a generalization of
the notion of split-graphs, \ref b.ST/. Call a $k$-uniform hypergraph $G$ 
a
{\it split-hypergraph} if its vertices can be partitioned into two 
classes
so that every $k$-tuple in one of the classes is an edge and no
$k$-tuple in the other class is an edge in $G$. We call such a
partition a splitting. Note that our definition allows one of the classes 
to have size smaller than $k$, in which case the requirement on that class 
automatically holds.

When $G$ is 
not a split hypergraph, we say that it 
is 
not split or non-split. 
For a hypergraph $G$, $V(G)$ will denote its vertex set and $E(G)$ its 
edge
set. 
From now on, $n$ denotes $|V(G)|$
and $N$ is $|V(G)|+|E(G)|$.

  The question arose naturally: can we state something similar to \ref 
t.graph/
for the case of split hypergraphs?
 {\it{U. N. Peled}} asked whether the class of $k$-uniform split
hypergraphs can be characterized by a finite list of excluded induced
subhypergraphs. We give an affirmative answer to this question in the 
next section.
Section 3 gives examples of minimal non-split 
hypergraphs in the 3-uniform case. We did not strive for the completeness of the
list. Nevertheless, it may show that the number of excluded induced 
subhypergraphs increases considerably as we go from graphs to 3-uniform 
hypergraphs.
In Section 4 we 
examine if the connection
between degree sequences and the split property of graphs can be 
generalized. We find a
negative answer. 

In their paper \ref b.ST/, where they apply split 
hypergraphs in 
learning
theory, 
{\it{R. H. Sloan}} and {\it{Gy. Tur\'an}} gave an algorithm that decides if a
$k$-uniform hypergraph is split and gives a splitting if possible. The 
running time of their
algorithm is
$O(n^{2k})$; they mention that finding a more
efficient algorithm is an open question. Feder, Hell, Klein and Motwani 
introduce the notion of a sparse and dense partition of a graph, which is 
a kind of generalization of a splitting. Their definition can be extended 
to hypergraphs and the algorithm they present (in Theorem 3.1) can be used 
to find a splitting (or all splittings) of a uniform hypergraph. The 
running time of their 
algorithm in this context is the same as the time requirement of the one 
given by 
Sloan and Tur\'an. In the last section we give an $O(N \log N)$ algorithm 
for the problem. Our method can also be used to give an algorithm for 
sparse 
and dense partitions that is faster than the one in \ref b.FHKM/ in 
certain cases.

The restriction of an edge to an $X\subset V(G)$ is its 
intersection with $X$. Minimality of a graph with respect to a certain 
property will 
always mean that no {\it induced} proper subgraph of it has the property.
As usual, $X \Delta Y$ will stand for the symmetric difference between  
sets $X$ and $Y$. When $x\in V(G)$, we sometimes simply use $x$ 
instead of 
$\{x\}$ in set operations.

\bsection{The existence of a finite characterization}{s.main}

In this section we prove our main theorem. First we state two lemmas. The 
following lemma holds for hypergraphs with edges of sizes {\it at most} 
$k$, which is a relaxation of $k$-uniformity, and is convenient for our 
inductive proof.

\procl l.mcs Let $k$ and $m$ be positive integers and $\Gamma$ a 
hypergraph 
with edges of 
size
at most $k$. Then $\Gamma$ has at most $k^m$ different minimal cover sets 
of at most
$m$ points. Furthermore, there is an algorithm of time requirement 
$O(|V(\Gamma)|+|E(\Gamma)|)$ that outputs a list of these minimal cover 
sets.
\endprocl

Because of the short running time of our algorithm, some more details 
about the implementation could be added. However, 
designing data 
structure and the
implementation of the algorithm does not match the texture of
our paper, so we leave these to the interested reader.

\proof We prove the first part of the claim by induction on $m$. The case 
$m=1$ is obvious.  Suppose that
for $m<r$ the statement is true, and let $\Gamma$ be any hypergraph
with edges of size at most $k$. Let $H$ be an arbitrary edge.  Any
\mlp of $\Gamma$ with at most $m$ points consists of a point $x$ of $H$
and a \mlp of $\Gamma |_{V(\Gamma )\backslash \{x\}}$ with at most $m-1$
elements. Applying the induction hypothesis to the latter one and noting
that $x$ could be chosen in $k$ different ways, we get the statement
for $m=r$. 

The algorithm will imitate the proof in the previous paragraph. For 
$i=1,\ldots,m$, pick the first edge $e_i$ that has not been deleted yet, a 
vertex $v_i$ of $e_i$, delete the 
edges that contain 
$v_i$, and let $i:=i+1$. The output is $\{v_1,\ldots, v_m\}$. Since 
the number of minimal cover sets is bounded by a constant, we can list all 
of them by systematically going through all possible ways of picking 
$v_i$ in 
the sequence of $e_i$, and still use 
only linear time. 
\Qed 

In what follows, for the sake of simplicity, the edges of the $k$-uniform 
\hg
$G$ will be referred to as {\it type 0 edges}, and the $k$-tuples not 
contained
in $E(G)$ as {\it type 1 edges}. (Thus we identify $G$ with a labeling 
of 
the 
edges of a complete
$k$-uniform \hg by numbers 0 and 1.)

Consider some minimal non-split hypergraph $G$. Minimality implies that 
for 
any $x\in G$ there is a splitting $(A_0^x,A_1^x)$ of $G\setminus x$, i.e. 
$A_0^x$ induces only edges of type 0, $A_1^x$ induces only edges of type 
1. (A class that only contains edges of type $i$ will sometimes be called a 
{\it type-i-class}.)

\procl l.close Fix a splitting $(A_0^v,A_1^v)$ of $G\setminus    
v$ for every $v\in V(G)$.
If $x,y\in V(G)$ are vertices in a minimal non-split 
hypergraph $G$ then $|A_i^x\cap A_{1-i}^y|<k$ and $|A_i^x\Delta 
A_i^y|\leq 
2k$ for $i=0,1$.\endprocl

\proof The first inequality follows by the fact that a $k$-tuple cannot 
be contained in the type $i$ class of a splitting (of $G\setminus x$) and 
a 
type $1-i$ class of a splitting (of $G\setminus y$) at the same time. For the 
second inequality, use the facts 
$$A_i^x\setminus A_i^y\subset A_i^x\cap(A_{1-i}^y\cup\{y\}),$$ 
and $|A_i^x\cap (A_{1-i}^y\cup\{y\})|\leq k$ by the first inequality; 
similarly with $x$ and $y$ interchanged.\Qed

In the proof we shall consider sets subindexed by 0 or 1. Such a set, of 
subindex $i$, will always be a class of edges of type $i$, or some 
``preliminary set" in constructing a class only with edges of type $i$. 
This notation will provide us with the following convenience in 
terminology. Given some set $S\subset V(G)$ by a name
subindexed by 
$i$ ($i=0,1$), let the {\it bad edges} of $S$ be the edges of type $1-i$ 
on $S$.

\procl t.main If a $k$-uniform non-split \hg has more than 
$4(k+1)k^{4k(k+1)+1}$ vertices then it has a proper 
induced subhypergraph that is non-split. 
\endprocl

Fix $k$. We shall give an algorithm. Its brief description is given in the 
next two paragraphs. The input is 
$G,x,(A_0^x,A_1^x), \alpha, \beta$, where $G$ is a minimal non-split 
hypergraph,
$x\in 
V(G)$,  
$(A_0^x,A_1^x)$ is a splitting of $G\setminus x$, 
$\alpha\in\{1,\ldots,4k\}$, $\beta\in\{1,\ldots,{k+1}\}$. The algorithm 
either 
stops without any output, or outputs 
a splitting of $G\setminus y$ for some $y\in V(G)$. It has a randomized 
part, 
to simplify description. 

At the beginning our algorithm adds $x$ to $A_0^x$, creating 
some type 1 edges (bad edges) in $B_0:=A_0^x\cup x$. There is another class 
$B_1$, set to be $A_1^x$ at the start. Then we choose a random \mcs of size at 
most $k+1$ for the bad edges, and move it over from $B_0$ to $B_1$. In 
general, at 
each step at most one of $B_0$ and $B_1$ contains any bad edges, depending on the 
parity of the step. The algorithm always chooses uniformly a random \mcs 
of size 
at most $k+1$ 
for these bad edges, and consisting of points that have not been moved  
between the two classes yet, and any such set is chosen with positive 
probability. (If there is no such set, then the algorithm fails.)
Move this \mcs to 
the other class. When we arrive at the $\alpha$'th step, let the 
$\beta$'th element of the \mcs be called $y$.
(We may assume that there is some fixed ordering on any subset of the 
vertices, 
for example the one given by the listing of the vertices in the input.)
Remove $y$ from the \mcs 
(and hence from both classes in the future), and continue the procedure. We 
run the algorithm for at most $4k$ steps. With positive probability, at 
some step, 
$(B_0,B_1)$ will be a splitting of $G\setminus y$, and this is the output.
 
The inputs $\alpha$ and $\beta$, together with the choices that the 
random part of the algorithm makes, determine $y$. Using \ref l.mcs/, we 
shall give a bound on the number of possible choices, and finally get the 
general
upper bound in \ref t.main/ for the number of vertices in $G$. 

Here is a formal presentation of the algorithm.
\bigskip

{\tt Algorithm 1 

FUNCTION BadEdges($i,B_i$) returns the list of type $1-i$ (bad)
edges
of $B_i$.

FUNCTION MinCoverSet($R,L$) returns a uniformly chosen random 
minimum
cover
set of $R$ of size at most $k+1$ that is disjoint from $L$. If 
there is no such set then the algorithm ends and outputs 
``Unsuccessful''.

Start of Algorithm

INPUT: a $k$-uniform hypergraph $G$, vertex $x$ of $G$, a 
splitting $(A_0^x,A_1^x)$ of $G\setminus x$, numbers 
$\alpha\in\{1,\ldots,4k\}$ and $\beta\in\{1,\ldots,{k+1}\}$.

$B_0:=A_0^x\cup\{x\}$, $B_1:=A_1^x$, $L:=\emptyset$, $i:=0$, $j:=1$, 
$\ell:=\emptyset$
 
 {\bf while} $j\leq 4k$ {\bf do} $\{$

  \quad $R:=$BadEdges($i,B_i$)

  \quad {\bf if} $R=\emptyset$ {\bf then} {\bf output}($B_0,B_1$), {\bf end} of 
algorithm 

  \quad $\ell :=$MinCoverSet$(R,L)$

  \quad $B_i:=B_i\setminus \ell$, $L:=L\cup\ell$

%
%
%
%
%
%
  
  \quad {\bf if} $j=\alpha$ {\bf then} $\{$ if 
$|\ell|<\beta$ {\bf then} {\bf output}(``Unsuccessful''), {\bf end} of 
algorithm, 

  \qquad {\bf else} $\ell:=\ell\setminus\{y\}$, where $y$ is the 
$\beta$'th element of $\ell$ $\}$
  
  \quad $B_{1-i}:=B_{1-i}\cup\ell$

  \quad $i:=1-i$, $j:=j+1$

  $\}$

{\bf output}(``Unsuccessful''), {\bf end} of 
algorithm}

\bigskip


\procl l.correctness
Fix $G,x$ and $(A_0^x,A_1^x)$ in the input.
For any $y\in V(G)$, {\tt Algorithm 1} outputs a splitting of 
$G\setminus y$ with 
positive 
probability for some values of $\alpha,\beta$.\endprocl

\proof Fix $y$, and a splitting $(A_0^y,A_1^y)$ of $G\setminus y$. None of the 
bad edges in $A_0^x\cup \{x\}$ can be contained in $A_0^y$, so there 
exists some
minimal cover set of these bad edges
that is contained in $A_1^y\cup\{y\}$. 
(Such a 
\mcs is chosen by {\tt MinCoverSet} with positive probability.) This 
minimal cover set 
cannot have more than $k+1$ elements, because then a $k$-tuple would be  
present in both $A_0^x$ and $A_1^y$, 
contradicting \ref l.close/. 
Since the points that we
moved are in $A_1^y\cup\{y\}$ by our choice, they get to their 
final class by the move,
so we may assume that they are not moved in later steps of the algorithm.
Hence, when the algorithm chooses a minimal 
cover set of the bad edges of size at most $k+1$ and moves this set from 
$A_0^x\cup\{x\}$ to 
$A_1^x$, (call the resulting sets 
$C_0$ and $C_1$ respectively), we have $|C_0\Delta A_0^y|<|A_0^x\Delta 
A_0^y|$, and 
also $|C_1\Delta A_1^y|\leq|A_1^x\Delta A_1^y|$. These inequalities are 
crucial in that we succeed after at most $4k$ repetitions.
When $y$ is in $\ell$, if $\alpha$ is equal to the actual $j$ and 
$y$ is the $\beta$'th element in $\ell$, then $y$ is removed from 
$\ell$, and its removal ensures that the remaining graph has a 
splitting. (This possible removal causes the second inequality above to be 
not necessarily strict.)

Similarly, denote by $D_0$ and $D_1$ the sets $B_0$ and $B_1$ respectively 
at some 
step, and let $C_0$ and $C_1$ be $B_0$ and $B_1$ respectively in the 
next step. Let 
$\ell$ be the set that was moved, as in the algorithm; we may 
assume by symmetry that it was moved from $D_0$ to $D_1$. So, 
$C_0=D_0\setminus \ell$ and $C_1 =D_1\cup(\ell\setminus y)$. There 
were 
some bad edges in $D_0$, because the algorithm did not stop. None of 
these bad 
edges can be contained in $A_0^y$, so there exists some
minimal cover set of the bad edges 
that is contained in $A_1^y$. 
We may assume that the elements of this minimal cover set are points that 
have 
not been moved at any previous step (i.e., they are not in $L$), because 
we are considering the outcome of the algorithm (of positive 
probability), when the 
moved vertices must reach the type of class that contains them in 
($A_0^y,A_1^y$) by the 
move. (The only vertex not in any of $A_0^y$ or $A_1^y$ is $y$, but it is 
removed from $\ell$ if $\alpha$ and $\beta$ are suitable.)
The algorithm chooses such a minimal cover set $\ell$ with positive 
probability, and then we have $$|C_0\Delta A_0^y|<|D_0\Delta A_0^y|\;\; 
{\rm 
and}\;\;|C_1\Delta A_1^y|\leq |D_1\Delta A_1^y|.$$ 

We conclude that $|B_0\Delta A_0^y|+|B_1\Delta A_1^y|$ decreases in 
each iteration step by at least 1 with positive probability. 
At some step $y$ has to be contained in $\ell$, since $G$ itself is not 
split. If $\alpha$ and $\beta$ are suitable, then $y$ is removed from $\ell$ 
(and the future $B_0\cup B_1$). Then the 
algorithm can 
stop for two reasons. The first one is when it outputs ``Unsuccessful'' 
because of the 
random choices that {\tt MinCoverSet} made. However, we have seen that for 
a certain sequence of the random choices that {\tt MinCoverSet} makes along the 
iteration steps, there always exists a minimal cover 
set of the bad edges consisting of points not moved yet, as long as the 
set of 
the
bad edges is nonempty. On the other hand, when the set of bad edges becomes 
empty, the actual $B_0$ and $B_1$ is a splitting. 
This either happens when the algorithm has already 
run for $4k$ steps ($j>4k$), 
in which case $|B_0\Delta A_0^y|+|B_1\Delta A_1^y|$ has to be 0 (at the first 
step $|B_0\Delta A_0^y|+|B_1\Delta
A_1^y|=|A_0^x\Delta A_0^y|+|A_1^x\Delta A_1^y|\leq 4k$,
by \ref l.close/), and we achieved the 
splitting $(A_0^y,A_1^y)$ as an output. Or there are no bad edges 
at some earlier step. In that case the actual $(B_0,B_1)$ does not necessarily 
coincide with the fixed splitting 
$(A_0^y,A_1^y)$, but it is still a splitting of $G\setminus y$.
\Qed

\proofof t.main Fix again $G$, $x$ and $(A_0^x,A_1^x)$.
According to \ref l.correctness/, for every $y\in V(G)$ there is 
some input $\alpha,\beta$ such that with positive probability {\tt 
Algorithm 1} gives a splitting of $G\setminus y$. The
corresponding sequence of $\ell$'s (along the iteration steps of the 
algorithm) and the $\alpha,\beta$ 
determine 
the output and $y$. There are at most 
$k^{k+1}$ choices for $\ell$ in each step, by 
\ref l.mcs/, and the rest of the algorithm is deterministic. So there are 
at most 
$k^{4k(k+1)}$ possible sequences of $\ell$'s 
along the iteration steps of the algorithm. There are $4k({k+1})$ 
possible 
inputs $\alpha,\beta$. 
Hence there are at most $k^{4k(k+1)}4k$ possible $y$'s, which is our 
upper bound 
for the number of vertices in a minimal non-split $k$-uniform 
hypergraph.\Qed

A straightforward corollary is our main theorem:

\procl t.th The family of $k$-uniform split hypergraphs can be 
characterized by a finite set of forbidden induced 
subhypergraphs.\endprocl

\bsection{3-uniform minimal non-split hypergraphs}{s.examples}

Although \ref t.main/ gave a 
very rough upper bound on the size of 
minimal 
non-split hypergraphs, in this section
we illustrate with a few examples that the list of excluded induced
subhypergraphs increases considerably compared to the case of graphs.

For $|V(G)| = 6$, $G$ is non-split if and only if for any edge its 
complement is also an edge
and if the subhypergraph induced by any 4 points is neither the empty
nor the complete 3-uniform hypergraph. It is clear that these are 
minimal 
non-split hypergraphs too.

We get examples of minimal non-split hypergraphs for the
case of $|V(G)| = 7,8,9$ in the following manner. Write the vertices
along a circle and let the edges be exactly those 3-tuples that consist of
three consecutive nodes along the circle. $G$ is not split, since we could put
into the type-0-class at most 3 vertices, but then the rest
contains a type 0 edge (and similarly if we put only 2 arbitrary points
in the type-0-class). However, dropping any vertex, two suitably chosen
points can cut up the remaining arc to parts with at most two vertices. In 
the case of 8 (respectively 9) vertices, we can erase 1 (respectively
1,2,3 or 4) edges from the \hg just described, and 
similarly, we
still get a minimal non-split hypergraph. The 
argument is nearly the same for these new graphs.

We have found minimal 3-uniform non-split hypergraphs on vertex sets of
11 and 12 elements, but the details are lengthy. The
examples above may show that even in the 3-uniform case there are more
than a hundred (and probably even much more) minimal excluded 
subhypergraphs.

\bsection{Degree sequences of hypergraphs}{s.degrees}

The question arises, whether a characterization like that in  part (iii) 
of
\ref t.graph/ can be given for split hypergraphs. Since the equivalence
of (iii) and (i) provides us with an $O(n)$ algorithm to decide if a 
graph
is split or not, a characterization for split
hypergraphs using degree sequences could be promising.

\procl p.degree There are 3-uniform hypergraphs $G$ and $G'$ on the vertex 
set $\{1,
\ldots, n\}$ such that $G$ is split, $G'$ is non-split, moreover
$d_i=d'_i$ and $\delta_{ij} = \delta'_{ij}$ for any $i,j\in
\{1,2,\ldots, n\}$. Here $d_i$ ($d'_i$) stands for the degree of
$i$ in $G$ ($G'$); $\delta_{ij}$ ($\delta'_{ij}$)
is the number of edges containing both $i$ and $j$ in $G$ ($G'$).
\endprocl

\proof
Call a 6-element subset of the vertices of a 3-uniform
\hg $H$ exchangeable if there is a list $a_1, a_2, a_3, a_{-1}, a_{-2},
a_{-3}$ of these vertices so that
$\{a_1, a_2, a_3\}$, $\{a_1, a_{-2}, a_{-3}\}$, $\{a_{-1}, a_2,
a_{-3}\}$, $\{a_{-1}, a_{-2}, a_3\}$ are edges of $H$, and the other 
3-tuples not
containing $a_i$ and $a_{-i}$ together for some $i$ are not in $H$. We say 
that we
exchange the exchangeable point set $\{a_1, a_2, a_3, a_{-1}, a_{-2}, 
a_{-3}\}$ when we
erase the four edges described above and add the 3-tuples $\{a_{-1}, 
a_{-2},
a_{-3}\}$, $\{a_{-1}, a_2, a_3\}$, $\{a_1, a_{-2}, a_3\}$, $\{a_1,
a_2, a_{-3}\}$ to $E(H)$.
We may think of the exchangeable 6 vertices as the nodes of
an octahedron. Then being exchangeable means that, if the nodes were 
labeled 
appropriately, the four 3-tuples
in $E(H)$ are determined by four faces whose pairwise
intersection is one point. Exchanging means that we replace these 
3-tuples in $E(H)$ by the other four faces of the octahedron.
One 
can easily see that these are well-defined, the listing of the 6 vertices
is essentially unique. On the other hand the degrees and the number of
edges incident to two points does not change after applying an
exchanging.  Thus if we get a non-split \hg from a split \hg after
serial exchanging, then the proposition follows. 

Now, denote
$X=\{1,2,3,4\}$, $Y=\{5,6,7,8\}$. Define $G$ on the vertex set
$\{1, \ldots, 14\}$ as follows. Every 3-tuple of $X\cup \{9,
10, 11\}$ is red, every 3-tuple of $Y\cup \{12,13,14\}$ is blue, and
if a 3-tuple intersects both $X$ and $\{12, 13, 14\}$ then it is 
blue, if a 3-tuple intersects both $Y$ and $\{9,10,11\}$ then it is red. 
Moreover, the 3-tuples of $\{9, 10, \ldots, 14\}$ that have not been 
defined
yet get colors so that $\{9, 10, \ldots, 14\}$ give an
exchangeable set of points, as in the first sentence of this proof 
with 
$(a_1,a_2,a_3,a_{-1},a_{-2},a_{-3})=(9,10,11,12,13,14)$. The colors of the 
remaining 3-tuples 
are
arbitrary.  $G$ is obviously split with color classes $X\cup \{9, 10,
11\}$ as red and $Y\cup \{12, 13, 14\}$ as blue. Let $G'$ be the \hg
obtained from $G$ by exchanging $\{9, \ldots, 14\}$. In $G'$,  $\{9, 10, 
11\}$ has become blue, so at least one of
its points must be put to the other class. But there this point gives a 
red
edge with any two points of $Y$. Since at least two points of $Y$
must remain in the blue class, we conclude that $G'$ is indeed non-split.  
\Qed

For an arbitrary $k$, one can give a $k$-uniform counterexample with a
construction similar to the one above, using the $k$-dimensional
cross-polytope $\{x\in\R^k\, :\, ||x||_1\leq 1\}$ (the ``octahedron" of 
dimension $k$).

\bsection{Algorithm for splitting}{s.algorithm}

\procl t.algorithm For an arbitrary, but fixed $k$ there is an algorithm 
that decides 
if
a $G$ $k$-uniform hypergraph is a split hypergraph, and gives a
splitting when it is. The running time is $O(N \log N)$ (where 
we defined $N=|V(G)| + 
|E(G)|$).
\endprocl

As in the case of \ref l.mcs/, we mention that the 
running time is understood with appropriately chosen data structure.
Let us also point out that the constant in the $O(N \log N)$ bound is                     
exponential in $k$. 

The following lemma is used repeatedly, to obtain splittings of larger and 
larger subhypergraphs of $G$.

\procl l.close2 Let $H$ be a $k$-uniform hypergraph, and $(V^1, V^2)$ be 
a partition of
$V(H)$. Denote by $H^1$ and $H^2$ the
subhypergraphs of $H$ induced by $V^1$ and $V^2$ respectively. 
Suppose that $(A_0,A_1)$ is a splitting of $H$ and 
$(A_0^i, A_1^i)$ 
is a splitting of $H^i$ $(i=1,2)$. Then $|A_j\backslash (A_j^1 \cup 
A_j^2)| \le 2k-2$, $j=0,1$.\endprocl

\proof Otherwise some $k$-tuple in $A_j$ would be present in 
$A_{1-j}^1$ or $A_{1-j}^2$ ($j\in \{0,1\}$), contradicting the fact that 
both are 
monochromatic of different types.  \Qed

We present an algorithm whose existence proves \ref t.algorithm/. Call it 
{\tt Algorithm 2}. We do not give the pidgin Pascal program this time, 
because the main structure of the algorithm is very simple, and the more 
particular elements are similar to {\tt Algorithm 1} (except that 
the random part is replaced by a deterministic choice, as we 
shall see). Let me start with an overview of the algorithm.
We use notation $\langle m\rangle$ for the modulo 2 
value of an integer $m$.

The crucial component of our algorithm is a subroutine called {\tt Tree}. 
Its 
input is $(H,X_0,X_1)$, where $H$ is a $k$-uniform hypergraph, and $(X_0,X_1)$ 
partitions $V(H)$. {\tt Tree} will create labels for 
the vertices of a rooted tree $T$ of depth $2k-2$, where each inner 
node has 
$k^{2k-2}$ children. For each vertex, we assume that there is some 
ordering on its children, so we can talk about the ``$m$'th child" (when 
$m\leq k^{2k-2}$). 
Similarly, we assume that for any set of $k$-tuples of $G$, there is some 
fixed ordering on the minimal cover sets of size at most $2k-2$ of these 
$k$-tuples, so 
we can talk 
about the ``$m$'th minimal cover set" (when $m$ is small enough). 
Note that such an ordering of the minimal cover 
sets can be easily implemented by an algorithm that is linear in 
the input size, using \ref l.mcs/. Finally, fix 
some ordering on the vertices of $T$, so that vertices at smaller 
depth precede vertices at larger depth. 

Say that the root {\it has depth 0}, its children 
{\it have depth 1}, etc. The root $r$ will have label $(X_0^r, 
X_1^r):=(X_0,X_1)$. If $X_0^r$ contains some edge of type 1 in $G$, let 
$s:=0$. Otherwise let $s:=1$. 

Suppose that $v$ is a vertex in $T$ whose parent is $u$, and $v$ is at 
depth
$g$. Assume that $v$ is the $m$'th 
child of $u$. Let $v$ get label $(\ell_v,X^v_0,X^v_1)$, where $\ell_v$ is 
the $m$'th minimal cover set of size at most $2k-2$ of the bad edges in 
$X_{\langle g+s\rangle}^u$, let
$X_{\langle g+s\rangle}^v:=X_{\langle g+s\rangle}^u\setminus \ell_v$ and 
$X_{\langle 1+g+s\rangle}^v=X_{\langle 1+g+s\rangle}^u\cup \ell_v$.  
If $m$ is such that there is no $m$'th minimal cover set of size at 
most $2k-2$ for 
the actual 
bad 
edges, then delete $v$ and all its offspring from $T$. If there are no bad 
edges in $X_{g+s}^u$, then {\tt Tree} outputs the first such 
$(X_0^u,X_1^u)$ (in 
the ordering of the vertices), which is a 
splitting of the input graph $H$. If there are bad edges for all the 
leaves, then {\tt Tree} outputs ``No", indicating $H$ is not split.

{\tt Algorithm 2} receives input $G$, a $k$-uniform hypergraph.
In the preparational step it partitions 
$V(G)$ to sets $A_1,\ldots, A_{\lceil n/{k-1}\rceil}$, where each $A_i$ 
has $k-1$ elements with 
the possible exception of $A_{\lceil n/{k-1}\rceil}$. Then 
we call the iteration part of the algorithm, with input $\bigl( 
(A_1,\emptyset),\ldots,(A_{\lceil n/{k-1}\rceil},\emptyset)\bigr)$. 

In each 
cycle of the iteration part of {\tt Algorithm 2} there is an incoming list 
$\bigl( 
(X_0^1,X_1^1),$ 
$(X_0^2,X_1^2),\ldots,  (X_0^m,X_1^m)\bigr)$, where $(X_0^i,X_1^i)$ is a 
splitting 
of $G|_{X_0^i\cup X_1^i}$, and the $(X_0^i\cup X_1^i)_{i=1}^m$ give 
a partition of $V(G)$. Now form pairs 
$(Y^i_0,Y^i_1):=(X_0^{2i-1}\cup X_0^{2i}, 
X_1^{2i-1}\cup 
X_1^{2i})$ as $i=1,\ldots, \lfloor m/2 \rfloor$; if $m$ is odd then define 
$(Y^{\lceil m/2 \rceil}_0,Y^{\lceil m/2 \rceil}_1):=(X_0^m,X_1^m)$. Now, call 
{\tt Tree} for each $i\in \{1,\ldots,{\lceil m/2 \rceil}\}$, with input 
$(Y_0^i,Y_1^i)$. If any of these returns ``No", then {\tt Algorithm 2} 
ends and outputs 
``$G$ has no splitting". Otherwise, if {\tt Tree} returns $(A_0^i,A_1^i)$ at the 
$i$'th call then the input for the next iteration step is $\bigl(
(A_0^1,A_1^1),
(A_0^2,A_1^2),\ldots,  (A_0^{\lceil m/2 \rceil},A_1^{\lceil m/2 
\rceil})\bigr)$. The last cycle of the iteration is when $m=2$, that 
is, when we only call {\tt Tree} once. Suppose that then
it returns 
the pair $(A_0^1,A_1^1)$. In this case {\tt Algorithm 2} outputs this 
pair, 
which is a splitting of $G$.

\procl l.correctalg {\tt Algorithm 2}
provides a splitting of $G$ if there 
is any, 
and answers ``$G$ has no splitting" if $G$ is non-split.\endprocl

\proof We shall prove that {\tt Tree} indeed decides if there is a 
splitting of the 
hypergraph that it receives while {\tt Algorithm 2} 
is running and 
outputs a splitting if there is one. 
(Note that the hypergraphs that {\tt Tree} receives in {\tt Algorithm 2} 
have a very specific form.)
Then it is clear that the iteration 
steps give splittings 
for families of bigger 
and bigger graphs and finally for $G$, if $G$ is split.
So, as in \ref l.close2/, let $H$ be a $k$-uniform hypergraph, and 
$(V^1,V^2)$ be a partition of $V(H)$. Denote by $H^1$ and $H^2$ the
subhypergraphs of $H$ induced by $V^1$ and $V^2$ respectively.
Suppose that
$(A_0^i, A_1^i)$
is a splitting of $H^i$ $(i=1,2)$. We need to show that if the input is 
$(H,A_0^1\cup A_0^2,A_1^1\cup A_1^2)$ then {\tt Tree} finds a splitting of 
$H$ if there is any, and answers ``No" if $H$ is not split. Note that the 
inputs that {\tt Tree} can get in {\tt Algorithm 2} indeed have this 
special form.

Suppose that $H$ has a splitting, fix one such $(A_0,A_1)$. The process of 
finding a splitting of $H$ corresponds to a walk on $T$, starting from the 
root and going in each step towards a child. Each step of the walk 
corresponds to finding a \mcs for the bad edges 
in the class being examined actually, and putting this set over to the other 
class (starting with the two classes from the input of {\tt Tree}). By the 
same argument as in \ref l.correctness/ (using the analogue 
\ref l.close2/ of \ref l.close/ now), we shall conclude that if there is 
a splitting, 
we get it in $\leq 2k-2$ steps, or equivalently, at the latest when arriving 
to a leaf of 
$T$.

To be more detailed, suppose that $(A_0^1\cup A_0^2,A_1^1\cup 
A_1^2)=:(B_0,B_1)$ is not a splitting yet (otherwise we are done, the 
algorithm outputs the splitting after making the label for the root). Then 
there are type 1 edges in $B_0$ or there are type 0 edges in $B_1$ (call any 
of these a {\it bad edge}). We may assume by symmetry that $B_0$ has bad 
edges. Then some \mcs of these bad edges is contained in $A_1$, and the root 
has a child $v$ such that $\ell_v$ is this \mcs. The iteration part of 
{\tt Algorithm 2} moves over $\ell_v$ from 
$B_0$ to 
$B_1$, hence creating $X_0^v$ and $X_1^v$. We have $|X_0^v\Delta 
A_0|+|X_1^v\Delta A_1|<|B_0\Delta 
A_0|+|B_1\Delta A_1|$, similarly to the displayed line in the proof of 
\ref l.correctness/. Note that after the 
first step of the walk on $T$ only one of $X_0^v$ and $X_1^v$ can have bad 
edges.

Proceed similarly: when the walk in $T$ is in $u$ and the bad edges are 
in $X_i^u$ ($i\in\{0,1\}$), there is a \mcs for these bad edges that is 
contained in $A_{1-i}$, and such that it is disjoint from 
$\bigcup \ell_w$, where $w$ ranges through the ancestors of $u$.
(This latter assumption can be made because in each step we moved points 
that got to their ``final class" in $(A_0,A_1)$ by this move.)
Some 
child 
$v$  
of $u$ is such that this \mcs is $\ell_v$, and for the arising 
$(X_0^v,X_1^v)$ we 
again have $$|X_0^v\Delta
A_0|+|X_1^v\Delta A_1|<|X_0^u\Delta
A_0|+|X_1^u\Delta A_1|.$$ 

Since for the root $r$ we have $|X_0^r\Delta
A_0|+|X_1^r\Delta A_1|\leq 2k-2$ by \ref l.close2/,
in at most $2k-2$ steps we necessarily arrive to a vertex where 
the two casses contain no bad edges. At this step {\tt Tree} outputs the two 
classes, 
which is a splitting of $H$.

If $H$ is not split, then it is clear that there are bad edges in some 
class for any partition of the 
vertex set, so {\tt Tree} will eventually output ``No". \Qed

\procl l.runningtime The algorithm runs in $O(N \log N)$ time.\endprocl

\proof First, {\tt Tree} runs in linear time, because the size of the 
underlying tree is a constant (determined by $k$), and the labels 
consisting of minimal cover sets can be constructed in a number of 
steps that is linear in the input of {\tt Tree}, using the method in \ref 
l.mcs/. 

In each iteration step each element of $V(G)\cup E(G)$ is present in at 
most one of the subhypergraphs that some {\tt Tree} subroutine receives, 
and then in each iteration cycle there is a linear number of extra steps
when the iteration cycle creates what it returns from the outputs of the 
{\tt Tree} subroutines.
This is $O(N)$, and there are $O(\log N)$ 
iteration steps, so we get the claim.  \Qed

\medbreak
\noindent {\bf Acknowledgements.}\enspace
I am grateful to P\'eter Hajnal for his support and suggestions.
I also thank Csaba Bir{\'o} for his comments on the manuscript.

\startbib
\bibitem[1]{BLS} A. Brandst{\"a}dt, V. B. Le, J. P.~Spinrad.
{\it Graph classes: a survey}.
SIAM Monographs on Discrete Mathematics and Applications.
Society for Industrial and Applied Mathematics (SIAM),
Philadelphia, 1999.
\bibitem[2]{FHKM} T. Feder, P. Hell, S. Klein, R. Motwani. List 
partitions. {\it SIAM J. of Discrete Math.}, 16: 449-478, 2003.
\bibitem[3]{GJ} M. R. Garey, D. S. Johnson. {\it Computers and 
intractability. A
guide to the theory of NP-completeness.} A Series of Books in the
Mathematical Sciences. W. H. Freeman and Co., San Francisco, Calif.,
1979.

\bibitem[4]{Go} M. C. Golumbic. {\it Algorithmic Graph Theory and Perfect 
Graphs.} Annals of Disc. Math., 57, Elsevier, 2004. 
\bibitem[5]{GyL} A. Gy\'arf\'as, J. Lehel. A Helly-type problem in trees. 
{\it Combinatorial theory and its applications}, II (Proc. Colloq., 
Balatonf\"ured, 1969), 571-584, North-Holland, Amsterdam, 1970.

\bibitem[6]{HS} P.L. Hammer, B. Simeone. The splittance of a graph. Univ. 
of Waterloo, 
Dept. of Combinatorics and Optimization, Res. Report CORR 77-39, 1977.

\bibitem[7]{K} R. M. Karp. Reducibility among combinatorial
problems. {\it Complexity of computer computations} (Proc. Sympos., IBM
Thomas J. Watson Res. Center, Yorktown Heights, N.Y., 1972),
85-103, Plenum, New York, 1972.

\bibitem[8]{NRSS} R. N. Naik, S. B. Rao, S. S. Shrikhande, N. M. Singhi.
Intersection graphs of $k$-uniform linear hypergraphs.
{\it European J. Combin.} 3, no. 2, 159-172, 1982.

\bibitem[9]{ST} R. H. Sloan, Gy. Tur\'an.  Learning from 
incomplete boundary queries using split graphs and
hypergraphs (extended abstract). {\it Computational learning theory}, 
Jerusalem, 1997, Lecture Notes in Comput. Sci., 38-50, Springer, 1997.
\bibitem[10] \'A. Tim\'ar. Split hypergraphs. Diploma thesis at the 
University of Szeged (in 
Hungarian), 2000.

\endbib  

\bibfile{\jobname}
\def\noop#1{\relax}
\input \jobname.bbl

\filbreak
\begingroup
\eightpoint\sc
\parindent=0pt\baselineskip=10pt

Department of Mathematics,
University of British Columbia,
121-1984 Mathematics Rd.,
Vancouver, BC V6T1Z1, Canada
\emailwww{timar[at]math.ubc.ca}{}
\htmlref{}{http://www.math.ubc.ca/$\sim$timar/}
\endgroup

\bye